\documentclass[12pt]{article}
\usepackage{epsfig,graphicx,color}
\usepackage{amsmath}
\usepackage{ulem}
\usepackage{latexsym}
\usepackage{amstext}
\usepackage{epsfig, graphicx}
\usepackage[round]{natbib}
\newcommand{\mb}[1]{ \mbox{\boldmath$#1$} }
\newcommand{\ds}{\displaystyle}
\newcommand{\beq}{\begin{eqnarray}}
\newcommand{\eeq}{\end{eqnarray}}
\newcommand{\beqq}{\begin{eqnarray*}}
\newcommand{\eeqq}{\end{eqnarray*}}
\newcommand{\p}{\partial}

\newcommand{\eps}{\varepsilon}
\newcommand{\x}{\mbox{\boldmath$x$}}

\newcommand{\Q}{\mbox{\boldmath$Q$}}

\newcommand{\y}{\mbox{\boldmath$y$}}

\newcommand{\w}{\mbox{\boldmath$w$}}

\newcommand{\n}{\mbox{\boldmath$n$}}

\font\bb=msbm10 at 12pt
\def\rR{\hbox{\bb R}}

\begin{document}
\pagestyle{plain}
\begin{center}
{\large \textbf{{Oscillatory survival probability and eigenvalues of the non-self adjoint Fokker-Planck operator }}}\\[5mm]
D. Holcman \footnote{Group of Applied Mathematics and Computational Biology,, Ecole
Normale Sup\'erieure, 46 rue d'Ulm 75005 Paris, France. This research is
supported by an ERC-starting-Grant.}, Z. Schuss  \footnote{Department of
Mathematics, Tel-Aviv University, Tel-Aviv 69978, Israel.}
\end{center}
\date{}
\begin{abstract}
We demonstrate the oscillatory decay of the survival probability of the stochastic dynamics
$d\x_\eps=\mb{a}(\x_\eps)\, dt +\sqrt{2\eps}\,\mb{b}(\x_\eps)\,d\w$, which is activated by small
noise over the boundary of the domain of attraction $D$ of a stable focus of the drift
$\mb{a}(\x)$. The boundary $\p D$ of the domain is an unstable limit cycle of $\mb{a}(\x)$. The
oscillations are explained by a singular perturbation expansion of the spectrum of the Dirichlet
problem for the non-self adjoint Fokker-Planck operator in $D$
\[L_\eps u(\x)=\,\eps\sum_{i,j=1}^2  \frac{\p ^2\left[ \sigma ^{i,j}\left(\x\right)
u(\x) \right]}{\p x^i\p x^j}-\sum_{i=1}^2\frac {\p \left[ a^i\left(\x\right)
u(\x)\right]} {\p x^i} =-\lambda_\eps u(\x),\]
with $\mb{\sigma}(\x)=\mb{b}(\x)\mb{b}^T(\x)$.
We calculate the leading-order asymptotic expansion of all eigenvalues $\lambda_\eps$ for small
$\eps$. The principal eigenvalue is known to decay exponentially fast as $\eps\to0$. We find that
for small $\eps$ the higher-order eigenvalues are given by
$\lambda_{m,n}=2n\omega_1+mi\omega_2+O(\eps)$ for
$n=1,2,\ldots,\,m=\pm1,\ldots$, where $\omega_1$ and $\omega_2$ are explicitly computed constants.
We also find the asymptotic structure of the eigenfunctions of $L_\eps$ and of its adjoint
$L^*_\eps$. We illustrate the oscillatory decay with a model of synaptic depression of neuronal
network in neurobiology.
\end{abstract}

\section{Introduction}
The stochastic dynamics in $\rR^d$
\begin{align}
d\x_\eps(t)=\mb{a}(\x_\eps(t))\, dt +\sqrt{2\eps}\,\mb{b}(\x_\eps(t))\, d\w(t),\label{SDEo}
\end{align}
where $\w(t)$ is Brownian motion, serves as a model for a variety of physical, chemical,
biological, and engineering diffusion processes. The case of an isotropic constant diffusion matrix
$\mb{b}(\x)$, e.g. $\mb{I}$, and a conservative drift field $\mb{a}(\x)$ that is a gradient of a
potential, is often the overdamped (Smoluchowski) limit of the Langevin equation. When the potential
forms a well the exit problem is to evaluate the probability density function of the first passage
time of the trajectories $\x_\eps(t)$ of (\ref{SDEo}) from any point in the well to its boundary and
to evaluate its functionals in the small-noise limit $\eps\to0$. This problem, which represents
thermal activation over a potential barrier, has been extensively studied in the past 70 years and
is well understood. However, in damped systems, such as the Langevin equation, the drift field is
not conservative. This is also the case of phase tracking and synchronization loops in RADAR and
communications theory and other important engineering applications \citep*[Sections 8.4, 8.5, and
Chapter 10]{DSP}, \citep*{OPT}. In these models the non-conservative drift field $\mb{a}(\x)$  may
have a stable focus with a domain of attraction $D$. The exit problem is then much more complicated
than in the conservative case. In some models of neuronal activity \citep*{Holcman2006}, the drift
field $\mb{a}(\x)$  has a stable focus with a domain of attraction $D$, whose boundary $\ D$ is an
unstable limit cycle of the drift (see Fig.\ref{f:Focus}). Experimental data and Brownian dynamics simulations of this model indicate oscillatory decay of the survival probability in this model, that needs to be resolved.
In the non-conservative cases the principal eigenvalue and eigenvector of the Fokker-Planck operator
corresponding to (\ref{SDEo}) are real while those of higher order are complex valued, which may
cause oscillations in the probability density function of the first passage time $\tau$. Although in
the small noise limit the principal eigenvalue $\lambda_0$ and the mean first passage time
$\bar\tau$ are related asymptotically by
\begin{align}
\lambda_0\sim\frac{1}{\bar\tau}\hspace{0.5em}\mbox{for}\
\eps\ll1,\label{lambdatau0}
\end{align}
and the stationary (and quasi-stationary) exit point density on $\p D$ are the normalized flux of
the  principal eigenfunction $u_0(\y)$ of the Fokker-Planck operator, higher order eigenvalues and
eigenfunctions can cause discernible oscillations in the survival probability of $\x_\eps(t)$ in
$D$. This, as well as other problems, raise the question of where is the spectrum of the
Fokker-Planck non-self-adjoint elliptic operator? and how it depends on the structure of the
dynamics such as the drift.

{The Dirichlet problem for elliptic operators of the form
\begin{align}
Lu(\x)=\eps\mb{\sigma}(\x)\nabla\cdot\nabla u(\x)+\mb{a}(\x)\cdot \nabla u(\x)
\label{elliptic}
\end{align}
in bounded domains with sufficiently regular boundaries is self-adjoint when $\mb{a}(\x)$ is a
gradient, e.g., when $\mb{a}(\x)=\mb{0}$. The eigenvalues of $L$ in this case were computed
explicitly for simple geometries, such as the sphere, cube, projective sphere, and other analytical
manifolds \citep*{Chavel}. The asymptotic behavior of high-order eigenvalues (for
$\mb{a}(\x)=\mb{0}$) is known from Weyl's theorem \citep*{Weyl}.} This is not
the case, however, for non self-adjoint operators. Krein-Rutman's theorem \citep*{Krein} asserts
that the principal eigenvalue is simple and positive. More recent attempts at characterizing the
spectrum can be found {\em i.a.} in \citep*{Trefethen}, \citep*{Davies}, and
\citep*{Sjostrand}. Stochastic approaches based on {the} large deviation principle are summarized in \citep*{Freidlin}.

In the case of the Fokker-Planck Dirichlet problem, it is a singularly perturbed non self-adjoint
operator and the reciprocal of the principal eigenvalue is asymptotically the mean
first passage time to the boundary of the domain of a diffusion process, which
can be evaluated  asymptotically in the small noise limit \citep*{book},
\citep*{DSP} (see early attempts in \citep*[and references therein]{Friedman}. This expansion
represents the result of nearly 50 years of collective effort to derive a refined asymptotic
expansion based on the WKB approximation and matched asymptotics theory. Not much, however, is known
about higher order eigenvalues.

{In the present paper we consider the noisy dynamics \eqref{SDEo} confined in a domain $D$, as shown in Figures \ref{f:Focus} and \ref{f:Focus-Holcman-Tsodyks}. We demonstrate that for small driving noise the decay of the survival probability of a random trajectory in $D$ is oscillatory, due to the complex eigenvalues of the non-self-adjoint Dirichlet problem (\ref{elliptic}) in $D$. More specifically,}
the drift field $\mb{a}(\x)$ is assumed to have a stable focus in $D$, whose boundary $\p D$ is an
unstable limit cycle of $\mb{a}(\x)$. To state the main results, we use the following notation: $s$
is arclength on $\p D=\{\x(s): 0\leq s<S\}$, measured clockwise, $\n(\x)$ is the unit outer normal
at $\x\in\p D$, $B(s)=|\mb{a}(\x(s))|$, and $\sigma(s)=\n(\x(s))^T\mb{\sigma}(\x(s))\n(\x(s))$. The function $\xi(s)$ is defined in (\ref{Berneq3}) below.

Our main result for higher order eigenvalues is the asymptotic expression
\begin{align}
\lambda_{m,n}=2n\omega_1 +m\omega_2 i+O(\eps),\quad n=1,\ldots,\ m=\pm1,\pm2,\ldots,\label{lmni}
\end{align}
where the frequencies $\omega_1$ and $\omega_2$ are defined as
\begin{align}
\omega_1=\frac{\omega_2}{2\pi}\int_0^S\frac{\sigma(s)\xi^2(s)}{B(s)}\,ds \hbox{ and }
\omega_2=\frac{2\pi}{\ds\int_0^S\frac{ds}{B(s)}},\quad
\end{align}
which is found by studying the boundary layer near the limit cycle, where the spectrum is hiding (see section 4). The leading order asymptotic expansion of the principal eigenvalue $\lambda_0$ for small $\eps$ is related to the MFPT by (\ref{lambdatau0}), whose asymptotic structure was found in \citep*{MatkowskySchuss1982} and \citep*{DSP}. Section 3 contains {a new refinement of the WKB analysis that is} used in section \ref{s:Apparent} to demonstrate the oscillations in the survival probability and in the exit density. This result resolves the origin of the non-Poissonian nature of many phenomena, such as the times neurons stay depolarized in population dynamics (see discussion).
\section{The survival probability and the eigenvalue problem}
The exit time distribution can be expressed in terms of the transition probability density function
(pdf) $p_\eps(\y,t\,|\,\x)$ of the trajectories $\x_\eps(t)$ from $\x\in D$ to $\y\in D$ in time $t$. The pdf
is the solution of the Fokker-Planck equation (FPE)
\begin{align*}
\frac{\p p_\eps(\y,t\,|\,\x)}{\p t}=&\,L_{\y} p(\y,t\,|\,\x)\hspace{0.5em}\mbox{for}\ \x,\y\in
D\\
p_\eps(\y,t\,|\,\x)=&\,0\hspace{0.5em}\mbox{for}\ \x\in\p D,\ \y\in D,\ t>0\\
p_\eps(\y,0\,|\,\x)=&\,\delta(\y-\x)\hspace{0.5em}\mbox{for}\ \x,\y\in D{,}
\end{align*}
where $\mb{\sigma}(\x)=\mb{b}(\x)\mb{b}^T(\x)$. The Fokker-Planck operator $L_{\y}$ is given by
\begin{align}
L_{\y}u(\y)=&\,\eps\sum_{i,j=1}^2  \frac{\p ^2\left[ \sigma ^{i,j}\left(\y\right)
u(\y) \right]}{\p y^i\p y^j}-\sum_{i=1}^2\frac {\p \left[ a^i\left(\y\right)
u(\y)\right]} {\p
y^i}
\end{align}
and its adjoint is defined by
\begin{align}
L_{\x}^*v(\x)=&\,\eps\sum_{i,j=1}^2  \sigma ^{i,j}\left(\x\right)\frac{\p ^2 v(\x)
}{\p x^i\p x^j}+\sum_{i=1}^2 a^i\left(\x\right)\frac {\p v(\x)} {\p
x^i}
\end{align}

The non-self-adjoint operators $L_{\y}$ and $L_{\x}$ with homogeneous Dirichlet boundary conditions have the
same eigenvalues $\lambda_{n,m}$, because the equations are real and the eigenfunctions $u_{n,m}(\y)$ of $L_{\y}$ and $v_{n,m}(\x)$ of $L_{\x}^*$ { are bases that are bi-orthonormal in the complex Hilbert space such that
\begin{align}
\int_{D} \bar v_{n,m}(\y) L_{\y} u_{n,m}(\y)\,d\y = \int_{D}\bar u_{n,m}(\y)  L_{\y}^* v_{n,m}(\y)\,d\y= \delta_{n,m}.
 \end{align}}
The solution of the FPE can be expanded as
\begin{align}
 p_\eps(\y,t\,|\,\x)=e^{- \lambda_0t}u_0(\y)v_0(\x)+
\sum_{n,m} e^{-\lambda_{n,m}t}u_{n,m}(\y)\bar v_{n,m}(\x),\label{pepsunif}
 \end{align}
where $\lambda_0$ is the real-valued principal eigenvalue and $u_0,v_0$ are the corresponding
positive eigenfunctions, {that is, solutions of $L_{\x}(u_0)=-\lambda_0 u_0$ and
$L_{\y}^* (v_0)=-\lambda_0 v_0$, respectively.} The conditional probability density
function of the exit point $\y\in\p D$ and the exit time $\tau$ is given by
\begin{align}
{\Pr}\left\{\x_\eps(\tau)=\y,\tau=t\,|\,\x_\eps(0)=\x\right\}=
\frac{\mb{J}(\y,t\,|\,\x)\cdot \mbox{\boldmath$\nu(y)$}}
{\ds{\oint\limits_{\p D}}\mb{J}(\y,t\,|\,\,\x)\cdot \mbox{\boldmath$\nu(y)$}\,dS_{\y}},\label{exdens17}
\end{align}
where the flux density vector is given by
\begin{align}
J^i(\y,t\,|\,\x)=&\,a^i(\y) p_\eps(\y,t\,|\,\x)-\eps\sum_{j=1}^d \frac{\p
\left[\sigma^{i,j}(\y) p_\eps(\y,t\,|\,\x)\right]}{\p y^j}\nonumber\\
=&\,-\eps\sum_{j=1}^d\sigma^{i,j}(\y) \left[e^{- \lambda_0t}\frac{\p u_0(\y)}{\p y^j}v_0(\x)+
\sum_{n,m} e^{-\lambda_{n,m}t}\frac{\p u_{n,m}(\y)}{\p y^j}\bar v_{n,m}(\x)\right].
\label{flux17}
\end{align}
Here $\mb{\nu}(\y)$ is the unit outer normal vector at the boundary point $\y$. Note that due to the
homogeneous Dirichlet boundary condition the undifferentiated terms drop from (\ref{flux17}).
Equation (\ref{exdens17}) can be understood as follows. The normal component of the
flux density vector at time $t$ at the point $\y\in\p D$ is the joint probability of trajectories to survive in $D$ by time $t$ and to be absorbed in a unit surface element $\y+dS_{\y}$ at time $t$. The denominator in (\ref{exdens17}) is the absorption flux in $\p D$ at this time. It follows
that the normalized flux is the conditional probability to survive up to time $t$ and be absorbed in the surface element $\y+dS_{\y}$ at time $t$.

The survival probability of $\x_\eps(t)$ in $D$, averaged with respect to a uniform initial distribution, is
given in terms of the transition probability density function $p_\eps(\y,t\,|\,\x)$ of the trajectories
$\x_\eps(t)$ as
\begin{align}
 {\Pr}_{\scriptsize\mbox{survival}}(t)=&\,\frac{1}{|D|}\int\limits_D\Pr\{\tau>t\,|\,\x\}\,d\x=
 \frac{1}{|D|}\int\limits_D\int\limits_Dp_\eps(\y,t\,|\,\x)\,d\y\,d\x\nonumber\\
 =&\,e^{- \lambda_0t}+
\sum_{n,m} \frac{e^{-
 \lambda_{n,m}t}}{|D|}\int\limits_Du_{n,m}(\y)\,d\y\int\limits_D\bar v_{n,m}(\x)\,d\x .
 \end{align}
The pdf of the escape time is given by
\begin{align}
\Pr\{\tau=t\}=-\frac{d}{dt}{\Pr}_{\scriptsize\mbox{survival}}(t)=&\,\lambda_0e^{- \lambda_0t}+
\sum_{n,m} \frac{\lambda_{n,m}e^{-
 \lambda_{n,m}t}}{|D|}\int\limits_Du_{n,m}(\y)\,d\y\int\limits_D\bar v_{n,m}(\x)\,d\x.\label{pdf}
\end{align}

\section{Asymptotic expansion of the principal eigenvalue}\label{s:first-ev}
{This section summarizes \citep*[Section 10.2.6]{DSP}, which presents the asymptotic method for the case of the principal eigenvalue and the associated eigenfunctions $u_0(\x)$ and $v_0(\x)$. This method is the basis for the construction of the asymptotic expansion of all higher order eigenvalues and eigenfunctions of the problem at hand. It is presented here for completeness.}
\subsection{The field $\mb{a}(\x)$}\label{ss:a}
The local geometry of $D$ near $\p D$ can be described as follows. We denote by
$\x'$ the orthogonal projection of a point $\x\in D$ near the boundary.  The
signed distance to the boundary
\begin{align*}
\rho(\x)=\left\{\begin{array}{ccl} -|\x-\x'|&\mbox{for}& \x\in D\\
|\x-\x'|&\mbox{for}& \x\not\in D\\
0&\mbox{for}& \x\in \p D,\end{array}\right.
\end{align*}
defines $\n(\x)=\nabla\rho(\x)$ as the unit outer normal at $\x\in\p D$. Similarly,
the arclength on the boundary,  measured counterclockwise from a given boundary
point to the point $\x'$, defines $s(\x)$ for $\x\in D$ near the boundary and
defines $\nabla s(\x)$ as the unit tangent vector at $\x\in\p D$. Thus the
transformation $\x\to (\rho,s)$, where $\rho=\rho(\x),\ s=s(\x)$, is a 1-1 smooth map
of a strip near the boundary onto the strip $|\rho|<\rho_0,\ 0\leq s\leq S$, where
$\rho_0>0$ and $S$ is the arclength of the boundary.
\begin{figure}[ht!]
\centering
\includegraphics[width=2.5in,height=2.5in]{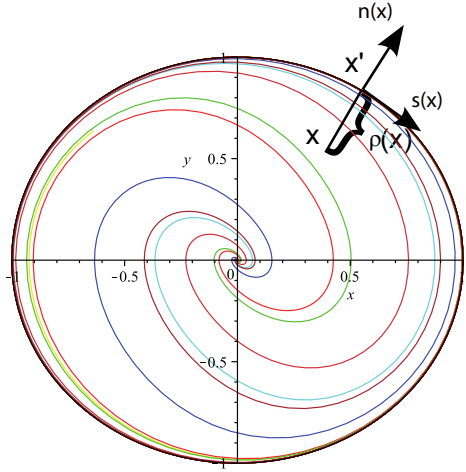}
\caption{\small The field $\mb{a}(\x)=[y,-x-y(1-x^2-y^2)]^T$ has a stable focus at the origin and
the boundary of the domain $D$ is a limit cycle.} \label{f:Focus}
\end{figure}
The transformation is given by $\x=\x'+\rho\nabla\rho(\x'),$ where $\x'$ is a function of $s$. The local representation of the field $\mb{a}(\x)$ in the boundary strip is assumed
\begin{align}
\mb{a}(\rho,s)=\left[a^0(s)\rho\nabla\rho+B(s)\nabla
s\right][1+o(1)]\hspace{0.5em}\mbox{for}\ \rho\to0,\label{arhos}
\end{align}
that is, the tangential component of the field at $\p D$ is
\beq
B(s)=\mb{a}(0,s)\cdot\nabla s=|\mb{a}(\x(s))|>0
\eeq
and the normal derivative of the normal component is $a^0(s)\geq0$ for all $0\leq s\leq S$. The decomposition (\ref{arhos}) for the field
$\mb{a}(\rho,s)$ in Figure \ref{f:Focus} is given by $a^0(s)=2\sin^2s,\ B(s)=1$.

\subsection{The { WKB} structure of the principal eigenfunction}\label{s:WKB1}
\subsubsection{The eikonal equation}
{We begin with the construction of the asymptotic approximation of the principal eigenfunction, now denoted  $u(\y)$. According to \citep*[Section 10.2.6]{DSP}, it has the WKB structure}
\begin{equation}
u(\y) = K_{\ds\eps}(\y) \exp\!\left\{-\frac{\psi(\y)}\eps \right\},
\label{WKBMD1}
\end{equation}
where the eikonal function $\psi(\y)$ is solution of the Hamilton-Jacobi {(eikonal)} equation
\begin{align}
\mb{\sigma}(\y)\nabla\psi(\y)\cdot\nabla\psi(\y)+\mb{a}(\y)\cdot\nabla\psi(\y)=0,\label{eikonalBF171}
\end{align}
{which is obtained by substituting (\ref{WKBMD1}) in (\ref{elliptic}) and comparing to
zero the leading term in the expansion of the resulting equation in powers of $\eps$
\citep*{MatkowskySchuss1977}, \citep*{book}, \citep*{MatkowskySchuss1982}. An interpretation of the
eikonal function $\psi(\y)$ in terms of the calculus of variations is given in large deviations
theory \citep*{Freidlin}.

The solution $\psi(\x)$ of the eikonal equation (\ref{eikonalBF171})} near the origin (the focus) is given by
\begin{align}
\psi(\x)=\frac12\x^T\Q\x+o(|\x|^2)\hspace{0.5em}\mbox{for}\ \x\to\mb{0}
\end{align}
with $\Q$ the solution of the Riccati equation
\begin{align}
2\mb{Q}\mb{\sigma}(\mb{0}) \mb{Q}+\mb{Q}\mb{A}+
\mb{A}^T\mb{Q}=\mb{0}.\label{Riccati}
\end{align}
where $\psi(\x)$ is the solution of the eikonal equation {(\ref{eikonalBF171})}.
The eikonal function $\psi(\x)$ is constant on $\p D$ with the local expansion
\begin{align}
\psi(\rho,s)=\hat\psi+\frac12\rho^2\phi(s)+o(\rho^2)\hspace{0.5em}\mbox{for}\
\rho\to0,\label{psirhos}
\end{align}
where $\phi(s)$ is the $S$-periodic solution of the Bernoulli equation
\begin{align}
\sigma(s)\phi^2(s)+a^0(s)\phi(s)+\frac12B(s)\phi'(s)=0\label{Berneq2}
\end{align}
and where $\sigma(s)=\mb{\sigma}(0,s)\nabla\rho(0,s)\cdot\nabla\rho(0,s).$  We may assume that for
isotropic diffusion $\sigma(s)=1$. Thus, for the dynamics in Fig.\ref{f:Focus}, the value of the constant $\hat\psi$ is calculated by integrating the characteristic equations for the eikonal equation (\ref{eikonalBF171}) \citep*{DSP}.

{To prove (\ref{psirhos}),} we note that $\psi(\y)$ is constant on the boundary, because in local coordinates on $\p D$ (\ref{eikonalBF171}) can be written  as
\begin{align}
[\nabla\psi(0,s)]^T\mb{\sigma}(0,s)\nabla\psi(0,s)+ B(s)\frac{\p \psi(0,s)}{\p
s}=0.\label{ee17}
\end{align}
To be well defined on the boundary, the function $\psi(0,s)$ must be periodic in
$s$ with period $S$. However, (\ref{ee17}) implies that the derivative $\p
\psi(0,s)/\p s$ does not change sign, because $B(s)>0$ and the matrix
$\mb{\sigma}(0,s)$ is positive definite. Thus we must have
\begin{align}
\psi(0,s)=\mbox{const.}=\hat\psi,\quad\nabla\psi(0,s)=0\hspace{0.5em}\mbox{for
all}\ 0\leq s\leq S.\label{psiconst}
\end{align}
It follows that near $\p   D$ the following expansion holds,
\begin{align*}
\psi(\rho,s)=\hat\psi+\frac12\rho^2\frac{\p ^2\psi(0,s)}{\p \rho^2}
+o\left(\rho^2\right)\hspace{0.5em}\mbox{as}\ \rho\to0.
\end{align*}
Setting $\phi(s)=\p ^2\psi(0,s)/\p \rho^2$ and using (\ref{arhos}) and
(\ref{psirhos}) in  (\ref{eikonalBF171}), we see that $\phi(s)$ must be the
$S$-periodic solution of the Bernoulli equation (\ref{Berneq2})
and $\sigma(s)=\mb{\sigma}(0,s)\nabla\rho(0,s)\cdot\nabla\rho(0,s)$. Writing $\xi_0(s)=\sqrt{-\phi(s)}$ in (\ref{Berneq2}),
we see that $\xi_0(s)$
is the $S$-periodic solution of the Bernoulli equation
\begin{align}
B(s)\xi_0'(s)+a^0(s)\xi_0(s)-\sigma(s)\xi_0^3(s)=0.\label{Berneq3}
\end{align}
{These function are discussed further in section \ref{ss:lambdatau}.}
\subsubsection{The transport equation}
The function $K_{\ds\eps}(\y)$ is a regular function of $\eps$ for $\y\in D$, but has to develop a boundary layer to satisfy the homogenous Dirichlet boundary condition
\begin{align}
K_{\ds\eps}(\y)=0\hspace{0.5em}\mbox{for}\ \y\in\p D.\label{KbcBF}
\end{align}
Therefore $K_{\ds\eps}(\y)$ is further decomposed into the product
\begin{align}
K_{\ds\eps}(\y)= \left[K_0(\y) + \eps K_1(\y) +\cdots\right]
q_{\ds\eps}(\y),\label{decomposeK}
\end{align}
where $K_0(\y),\, K_1(\y),\,\ldots$ are regular functions in $ D$ and on its
boundary and are independent of $\eps$, and $q_{\ds\eps}(\y)$ is a boundary
layer function. { As in the case of the eikonal equation $\psi(\y)$, the functions
$K_j(\y)\ (j=0,1,\dots)$ are solutions of first-order linear transport equations derived by
substituting (\ref{WKBMD1}) in (\ref{elliptic}), expanding the resulting equation in powers of
$\eps$, and equating to zero their coefficients \citep*{MatkowskySchuss1977}, \citep*{book}. These
functions
cannot satisfy the boundary condition (\ref{KbcBF}), because they are solutions of first-order
equations.} Thus $K_0(\y)$ has to be found by integrating {a transport equation along characteristics. Consequently, a boundary layer function $q_\eps(\y)$ is needed to make (\ref{decomposeK}) satisfy the homogeneous Dirichlet boundary condition.

The boundary layer function $q_{\ds\eps}(\y)$ has to satisfy} the boundary condition
\begin{align}
q_{\eps}(\y)=0\hspace{0.5em}\mbox{for}\ \y\in\p D,\label{qbcBF}
\end{align}
the matching condition
\begin{align}
\lim_{\eps\to0}q_{\eps}(\y) =&\,1\hspace{0.5em}\mbox{for all}\ \y\in D,
\label{qmatching}
\end{align}
and the smoothness condition
\begin{align}
\lim_{\eps\to0}\frac{\p ^iq_{\eps}(\y)}{\p (y^j)^i} =0,\hspace{0.5em}\mbox{for
all}\ \y\in D,\ i\geq1,\,1 \leq j\leq 2.\label{qsmoothness}
\end{align}
{First, we derive the transport equation for the leading term $K_0(\y)$.} The function
$K_{\ds\eps}(\y)$, which satisfies the transport equation
\begin{align}
&\,\eps \sum_{i,j=1}^2 \frac{\p ^2\sigma^{i,j}(\y) K_{\ds\eps}(\y)} {\p
y^i\p  y^j}\nonumber\\
&\,-\sum_{i=1}^2\left[2\sum_{j=1}^d\sigma^{i,j}(\y) \frac{\p \psi(\y)}{\p y^j}+
a^i(\y)\right] \frac{\p  K_{\ds\eps}(\y)} {\p  y^i}  \nonumber \\
&\,-\sum_{i=1}^2\left\{\frac{\p  a^i(\y)}{\p  y^i}+
\sum_{j=1}^2\left[\sigma^{i,j}(\y) \frac{\p ^2\psi(\y)} {\p  y^i\p
y^j}+2\frac{\p \sigma^{i,j}(\y)} {\p  y^j} \frac{\p \psi(\y)}{\p
y^j}\right]\right\}K_{\ds\eps}(\y)=0,\label{transportBF}
\end{align}
cannot have an internal layer at the global attractor point $\mb{0}$ in $ D$,
because stretching $\y=\sqrt{\eps}\mb{\xi}$ and taking the limit $\eps\to0$
converts the transport equation (\ref{transportBF}) to
\begin{align*}
&\,\sum_{i,j=1}^d \frac{\p ^2\sigma^{i,j}(\mb{0}) K_{0}(\mb{\xi})} {\p \xi^i\p
\xi^j} -(2\mb{A}\mb{Q}+\mb{A})
\mb{\xi}\cdot\nabla_{\mb{\xi}}K_{0} (\mb{\xi})\\
&\,-\mbox{tr}\left(\mb{A}+\mb{\sigma}(\mb{0})
 \mb{Q}\right)K_{0}(\mb{\xi})=0,
 \end{align*}
whose bounded solution is $K_{0}(\y)=const,$ because
$\mbox{tr}\left(\mb{A}+\mb{\sigma}(\mb{0})\mb{Q}\right)=0.$ The last equality
follows from the Riccati equation (\ref{Riccati}) (left multiply by
$\mb{Q}^{-1}$ and take the trace).

In view of eqs. (\ref{decomposeK})--(\ref{qsmoothness}), we obtain in the limit
$\eps\to0$ the transport equation\index{transport equation}
\begin{align}
&\,\sum_{i=1}^d\left[2\sum_{j=1}^d\sigma^{i,j}(\y) \frac{\p \psi(\y)}{\p  y^j}+
a^i(\y)\right] \frac{\p  {K_{0}(\y)}}
{\p  y^i}\label{transportsp} \\
 =&\,-\sum_{i=1}^d\left\{\frac{a^i(\y)}{\p  y^i}+
\sum_{j=1}^d\left[\sigma^{i,j}(\y) \frac{\p ^2\psi(\y)} {\p  y^i\p  y^j}+2%
\frac{\p \sigma^{i,j}(\y)} {\p  y^j}\frac{\p \psi(\y)}{\p
y^j}\right]\right\}{K_{0}(\y)}.  \nonumber
\end{align}
Because the characteristics diverge, the initial value on each characteristic
of the eikonal equation (\ref{eikonalBF171}) is given at $\y=\mb{0}$ as
$K_{0}(\mb{0})=const.$ (e.g., $const.=1$).

Note that using (\ref{arhos}) and (\ref{eikonalBF171}), the
field in the transport equation (\ref{transportBF}) can be written in local
coordinates near the boundary as
\begin{align}
2\mb{\sigma}(\y)\nabla\psi(\y) + \mb{a}(\y)=&\,2\mb{\sigma}(0,s)\nabla
\psi(0,s)+\mb{a}(0,s)+o(\rho)\nonumber\\
=&\,\rho\left[2\phi(s)\mb{\sigma}(0,s)\nabla\rho(0,s)+a^0(s)\nabla\rho(0,s)\right]
+o(\rho)\label{atilde17}
\end{align}
and the transport equation for $K_0(\y)$ can be written on $\p D$ as the linear
equation (which corrects eq.(10.125) in \citep*{DSP})
\begin{align}
B(s)\frac{dK_0(0,s)}{ds}+[a^0(s)+\sigma(s)\phi(s)+B'(s)]K_0(0,s)=0.
\end{align}
Using the relations (\ref{xiphi}) below, we obtain the solution
\begin{align}
K_0(0,s)=K_0\frac{\sqrt{-\phi(s)}}{B(s)},\label{K0s}
\end{align}
where $K_0=const.$ (e.g., $K_0=1$).
\subsubsection{The boundary layer equation for {$q_{\ds\eps}(\x)$}}\label{BLEQ}
To derive the boundary layer equation, we introduce the stretched variable
$\zeta=\rho/\sqrt{\eps}$ and define $q_{\ds\eps}(\x)= Q(\zeta,s,\eps)$.
Expanding all functions in (\ref{WKBMD1}) in powers of $\eps$ and
\begin{align}
Q(\zeta,s,\eps)\sim Q^0(\zeta,s)+\sqrt{\eps}Q^1(\zeta,s)+\cdots,\label{Qexp17}
\end{align}
and using (\ref{atilde17}), we obtain the boundary layer equation
\begin{align}
\sigma(s)\frac{\p ^2Q^0(\zeta,s)}{\p \zeta^2}-
\zeta\left[a^0(s)+2\sigma(s)\phi(s)\right] \frac{\p Q^0(\zeta,s)}{\p
\zeta}-B(s)\frac{\p  Q^0(\zeta,s)}{\p  s} =0.\label{ble17}
\end{align}
The boundary and matching conditions (\ref{qbcBF}), (\ref{qmatching}) imply
that
\begin{align}
Q^0(0,s)=0,\quad\lim_{\zeta\to-\infty} Q^0(\zeta,s )=1.\label{bcmathing17}
\end{align}
To solve (\ref{ble17}), (\ref{bcmathing17}), we set $\eta=\xi(s)\zeta$,
$Q^0(\zeta,s)=\tilde Q^0(\eta,s)$, and rewrite (\ref{ble17}) as
\begin{align}
\sigma(s)\xi^2(s)\frac{\p ^2\tilde Q^0(\eta,s)}{\p \eta^2}-&
\eta\left[a^0(s)+2\sigma(s)\phi(s)+\frac{B(s)\xi'(s)}{\xi(s)}\right] \frac{\p
\tilde Q^0(\eta,s)}{\p \eta}\nonumber\\
-&B(s)\frac{\p \tilde Q^0(\eta,s)}{\p  s} =0.\label{ble18}
\end{align}
Choosing $\xi(s)$ to be the $S$-periodic solution of the Bernoulli equation
(\ref{Berneq2}) the boundary value and matching problem (\ref{ble17}),
(\ref{bcmathing17}) becomes
\begin{align}
\frac{\p ^2\tilde Q^0(\eta,s)}{\p \eta^2}+ \eta \frac{\p\tilde Q^0(\eta,s)}{\p
\eta}-\frac{B(s)}{\sigma(s)\xi^2(s)}\frac{\p  \tilde Q^0(\eta,s)}{\p  s}
&\,=0,\label{ble18p}\\
\tilde Q^0(0,s)=0,\quad\lim_{\eta\to-\infty}\tilde Q^0(\eta,s
)&\,=1,\label{bcmathing17b}
\end{align}
which has the $s$-independent solution
\begin{align}
\tilde Q^0(\eta,s)=-\sqrt{\frac2\pi}\int\limits_0^{\eta}e^{-z^2/2}\,dz,\label{Q0int18}
\end{align}
that is,
\begin{align}
Q^0(\zeta,s)=-\sqrt{\frac2\pi}\int\limits_0^{\xi(s)\zeta}e^{-z^2/2}\,dz.\label{Q0int17}
\end{align}
{The uniform expansion of the first eigenfunction is constructed by putting together}
(\ref{WKBMD1}),  (\ref{decomposeK}), (\ref{Qexp17}), and (\ref{Q0int17}) {to} obtain that
\begin{align}
u_0(\y)=\exp\left\{-\frac{\psi(\y)}{\eps}\right\}\left[K_0(\y)+O(\sqrt{\eps})\right]
\sqrt{\frac2\pi}\int\limits_0^{\ds\frac{-\rho(\y)\xi(s(\y))}{\sqrt{\eps}}}e^{-z^2/2}\,dz,
\label{punif1752}
\end{align}
where $O(\sqrt{\eps})$ is uniform in $\y\in\bar  D$.

{Because $\eta=\xi(s)\zeta=\xi(s)\rho/\sqrt{\eps}$, equations (\ref{psirhos}) and
(\ref{xiphi}) near the boundary give}
\begin{align}
\psi(\rho,s)=\hat\psi-\frac{\rho^2\xi^2(s)}{2}+o(\rho^2),
\end{align}
{so} the eigenfunction (\ref{punif1752})  near the limit cycle has the form
\begin{align}
u_0(\y)\sim {\exp\left\{-\frac{\hat\psi}{\eps}\right\}}
\exp\left\{-\frac{\eta^2}{2}\right\}\left[K_0(\y)+O(\sqrt{\eps})\right]
\sqrt{\frac2\pi}\int\limits_0^{-\eta}e^{-z^2/2}\,dz \label{pbl}.
\end{align}
{The function $u_0(\y)$ is defined up to a multiplicative constant.}

The eigenfunction expansion (\ref{pepsunif}) and the expansion (\ref{punif1752}) of the principal eigenfunctions of the operator and its adjoint, respectively, give {the probability flux density}
\begin{align}
\mb{J}\cdot\mb{\nu}|_{\p D}(s,t)
\sim e^{-\lambda_0t}\sqrt{\frac{2\eps}{\pi}}
K_0(0,s)\xi(s)\sigma(s)e^{-\hat\psi/\eps}+\ldots,\label{Jdnu17}
\end{align}
hence, for $\y\in\p D${, which corresponds to $\rho=0$ and arclength $s$,}
\begin{align}
&\,\Pr\{\x(\tau)=\y,\tau=t\,|\,\x(0)=\x\}\nonumber\\
=&\,\frac{K_0(0,s)\xi(s)\sigma(s)+e^{(\lambda_0-\lambda_{n,m})t}u_{n,m}(\y)v_{n,m}(\x)+\ldots}{\ds\int_0^S
K_0(0,s)\xi(s)\sigma(s) \,ds+e^{(\lambda_0-\lambda_{n,m})t}u_{n,m}(\y)v_{n,m}(\x)+\ldots\,\ds}.\label{exitchat172}
\end{align}
Using (\ref{pbl}) at $\eta=0$ and (\ref{K0s}), we recover in the limit $t\to\infty$ the exit density
at $\y=(0,s)$ as \citep*{DSP}
\begin{align}
\Pr\{\x(\tau)=\y\,|\,\x\}\sim\frac{\ds\frac{\xi^2(s)\sigma(s)}{B(s)}}{
\ds{\int_0^S}\ds\frac{\xi^2(s)\sigma(s)}{B(s)}\,ds},\label{exitchat1721}
\end{align}
{which to leading order is independent of $\x$ outside a boundary layer of width $\sqrt{\eps}$.}
\subsubsection{The first eigenfunction of the adjoint problem}
The first eigenfunction $v_0(\x)$ of the backward operator $L^*_{\x}$ does not have the WKB
structure (\ref{WKBMD1}), but rather converges to a constant as $\eps\to0$, at every $\x\in D$
outside the boundary layer. Thus it {is merely the boundary layer $q_{\ds\eps}(\y)$}. Expanding as in section \ref{BLEQ}, we obtain the boundary value and matching problem
\begin{align}
\sigma(s)\frac{\p ^2Q^0(\zeta,s)}{\p \zeta^2}+ \zeta a^0(s) \frac{\p
Q^0(\zeta,s)}{\p \zeta}+B(s)\frac{\p Q^0(\zeta,s)}{\p  s} =0\label{ble19}\\
Q^0(0,s)=0,\quad\lim_{\zeta\to-\infty}Q^0(\zeta,s)=1.\label{BMC2}
\end{align}
The scaling $\eta=\xi_0(s)\eta$, with $\xi_0(s)$ the solution of
(\ref{Berneq3}), converts (\ref{ble19}), (\ref{BMC2}) to
\begin{align}
\frac{\p ^2\tilde Q^0(\eta,s)}{\p \eta^2}+ \eta \frac{\p\tilde Q^0(\eta,s)}{\p
\eta}+\frac{B(s)}{\sigma(s)\xi^2(s)}\frac{\p  \tilde Q^0(\eta,s)}{\p  s}
&\,=0,\label{ble20}\\
\tilde Q^0(0,s)=0,\quad\lim_{\eta\to-\infty}\tilde Q^0(\eta,s
)&\,=1,\label{bcmathing19}
\end{align}
where $\tilde Q^0(\eta,s)=Q^0(\zeta,s)$.
Using the solution of the Bernoulli equation (\ref{xiphi}), we obtain the $s$-independent solution (\ref{Q0int18}) and hence
(\ref{Q0int17}), which is the uniform approximation to $v_0(\y)$. We conclude that
\begin{align}
v_0(\y)={C_{\eps}}\mbox{erf}\left(\frac{\rho(\y)\xi(s(\y))}{\sqrt{\eps}}\right),
\label{punif1752b}
\end{align}
where $C_{\eps}$ depends on the normalization. Thus
\begin{align}
u_0(\y)\sim\exp\left\{-\frac{\psi(\y)}{\eps}\right\}\left[K_0(\y)+O(\sqrt{\eps})\right]v_0(\y),\label{u0}
\end{align}
which in the boundary layer coordinates has the form
\begin{align}
u_0(\y)\sim\exp\left\{-\frac{\eta^2}{2}\right\}\left[K_0(\y)+O(\sqrt{\eps})\right]v_0(\y).\label{v0}
\end{align}
\subsection{The principal eigenvalue $\lambda_0$ and the mean first passage time}
\label{ss:lambdatau}
{The asymptotic expansion of the mean first passage time $\bar\tau(\x)$ from $\x\in D$ to the boundary  is the solution of the Pontryagin-Andronov-Vitt boundary value problem
\citep*{MatkowskySchuss1982}, \citep*[Section 10.2.8]{DSP}}
\begin{align}
L^*\bar\tau(\x)=&\,-1\hspace{0.5em}\mbox{for}\ \x\in D\\
\bar\tau(\x)=&\,0\hspace{0.5em}\mbox{for}\ \x\in \p D.
\end{align}
{It is known to be independent of $\x$ outside the boundary layer in the sense that
\begin{align*}
\lim\limits_{\eps\to0}\frac{\bar\tau(\x)}{\bar\tau(\mb{0})}=1,
 \end{align*}
where}
\begin{align}
\bar\tau(\mb{0})\sim\frac{\pi^{3/2}\sqrt{2\eps\,\mbox{det}
\Q}}{\ds\int_0^{S}K_0(s)\xi(s)\,ds}
\exp\left\{\frac{\hat\psi}{\eps}\right\}.\label{tau}
\end{align}
The function $K_0(s)$, given by,
\begin{align}
K_0(s)=\frac{1}{B(s)}\exp\left\{-\int_0^s\left[\frac{a_0(s')-\xi^2(s')}{B(s')}\,ds'\right]\right\},
\end{align}
where $\xi(s)$ is the $S$-periodic solution of the Bernoulli equation
\begin{align}
\sigma(s)\xi^3(s)+[a^0(s)+2\sigma(s)\phi(s)]\xi(s)+B(s)\xi'(s)=0,\label{Berneqxi}
\end{align}
is defined up to a multiplicative constant that can be chosen to be 1. The solutions of the three Bernoulli equations $\phi(s)$ of (\ref{Berneq2}), $\xi(s)$ of (\ref{Berneqxi}), and $\xi_0(s)$
of (\ref{Berneq3})
are related to each other as follows (see \citep*[Section 10.2.6]{DSP} and Section 10.2.8),
\begin{align}
\xi_0(s)=\sqrt{-\phi(s)}=\xi(s).\label{xiphi}
\end{align}
{The mean first passage time from $\x\in D$ to the boundary is also given by
\citep*{DSP}}
\beq
\bar \tau(\x)&=& \int\limits_{0}^{\infty}
t\Pr\{\tau=t\,|\,\x\}dt=\int\limits_{0}^{\infty}\Pr\{\tau>t\,|\,\x\}= \int\limits_{0}^{\infty}
\int_{D} p_\eps(\y,t\,|\,\x)dt d\y \nonumber\\
&=&\frac1{\lambda_0}v_0(\x) + \sum_{n,m} \frac{\bar
v_{n,m}(\x)}{\lambda_{n,m}}\int\limits_Du_{n,m}(\y)\,d\y .\label{pdffpt}
\eeq
{If $\x$ is outside the boundary layer, then $v_0(\x)\sim 1$, as shown above, and $\int_Du_{n,m}(\y)\,d\y\sim 0$ by bi-orthogonality. Therefore
\beq
\bar \tau(\x) = \frac1{\lambda_0}v_0(\x) (1+o(1))\hspace{0.5em}\mbox{for}\ \eps\ll1.\label{pdffptbis}
\eeq
}
The principal eigenvalue $\lambda_0$ introduced in (\ref{lambdatau0}) is thus given more precisely by the asymptotic relation {
\begin{align}
\lambda_0\sim\frac{1}{\bar\tau(\mb{0})}\hspace{0.5em}\mbox{for}\
\eps\ll1\label{lambdatau}
\end{align}
and in view of (\ref{tau}), $\lambda_0$ decreases exponentially fast as $\eps\to0$.}
\section{Higher order eigenvalues}\label{s:HO-ev}
{The asymptotic expansion of higher-order eigenfunctions is constructed by the method
used above to derive that of the principal eigenfunctions. First, we consider higher-order eigenfunctions of the adjoint problem, which leads to the boundary layer equation and matching conditions}
\begin{align}
\frac{\p ^2\tilde Q^0(\eta,s)}{\p \eta^2}+ \eta \frac{\p\tilde Q^0(\eta,s)}{\p
\eta}+\frac{B(s)}{\sigma(s)\xi^2(s)}\frac{\p  \tilde Q^0(\eta,s)}{\p  s}
&\,=-\frac{\lambda}{\sigma(s)\xi^2(s)} \tilde Q^0(\eta,s),\label{bleHO}\\
\tilde Q^0(0,s)=0,\quad\lim_{\eta\to-\infty}\tilde Q^0(\eta,s
)&\,=0.\label{bcmathingHO}
\end{align}
Separating $\tilde Q^0(\eta,s)=R(\eta)T(s)$, we obtain for the even function $R(\eta)$
the eigenvalue problem
\begin{align}
R''(\eta)+\eta R'(\eta)+\mu R(\eta)=0,\quad
R(0)=0,\quad\lim_{\eta\to-\infty}R(\eta)=0,\label{blevp}
\end{align}
where $\mu$ is the separation constant. The large $\eta$ asymptotics of
$R(\eta)$ is $R(\eta)\sim\exp\{-\eta^2/2\}$, so the substitution
$R(\eta)=\exp\{-\eta^2/4\}W(\eta)$ converts (\ref{blevp}) to the parabolic
cylinder function eigenvalue problem
\begin{align}
W''(\eta)+\left(\mu-\frac12-\frac{\eta^2}{4}\right)W(\eta)=0,\quad
W(0)=0,\quad\lim_{\eta\to-\infty}W(\eta)=0.\label{blevpW}
\end{align}
The eigenvalues of the problem (\ref{blevpW}) are $\mu_n=2n,\ (n=1,2,\ldots)$ with the
eigenfunctions
$$W_{2n+1}(\eta)=\exp\left\{-\frac{\eta^2}{4}\right\}H_{2n+1}\left(\frac{\eta}{\sqrt{2}}\right),$$
where $H_{2n+1}(x)$ are the  Hermite polynomials of odd orders \citep*{abramowitz}. Thus the radial
eigenfunctions are
\begin{align}
R_n(\eta)=\exp\left\{-\frac{\eta^2}{4}\right\}W_{2n+1}(\eta)
=\exp\left\{-\frac{\eta^2}{2}\right\}H_{2n+1}\left(\frac{\eta}{\sqrt{2}}\right).\label{Rn}
\end{align}
{The associated function $T(s)$ (normalized to one) is the $S$-periodic solution of
\begin{align}
-\mu_n T(s)+\frac{B(s)}{\sigma(s)\xi^2(s)}\frac{\p  T(s)}{\p  s}&\,=-\frac{\lambda}{\sigma(s)\xi^2(s)} T(s),\label{bleHOT}
\end{align}
given by
\begin{align}
T_{n}(s)=\exp\left\{-\lambda\int_0^s\frac{ds'}{B(s')}+ 2n\int_0^s\frac{\sigma(s')\xi^2(s')}{B(s')}\,ds'\right\}.
\end{align}
We introduce therefore}  the period and angular frequency of rotation of
the drift about the boundary which are, respectively,
$${\cal T}=\int_0^S\frac{ds'}{B(s')},\quad\omega=\frac{2\pi}{\cal T}.$$
{Thus $S$-}periodicity implies the relation
\begin{align}
-\lambda\int_0^S\frac{ds}{B(s)}+2n\int_0^S\frac{\sigma(s)\xi^2(s)}{B(s)}\,ds=2\pi
mi
\end{align}
for $m=\pm1,\pm2,\ldots\,.$ {It follows that for $n=1,\ldots$} the eigenvalues are
\begin{align}
\lambda_{m,n}=\left[\frac{n}{\pi}\int_0^S\frac{\sigma(s)\xi^2(s)}{B(s)}\,ds+mi\right]\omega
\label{moi}
\end{align}
and the rotational eigenfunctions are
\begin{align}
T_{m,n}(s)=\exp\left\{-\lambda_{m,n}\int_0^s\frac{ds'}{B(s')}+
2n\int_0^s\frac{\sigma(s')\xi^2(s')}{B(s')}\,ds'\right\}.
\end{align}
The eigenfunctions $\tilde Q_{m,n}(\eta,s)=R_n(\eta)T_{m,n}(s)$ are given by
\begin{align}
\tilde Q_{m,n}(\eta,s)=\exp\left\{-\frac{\eta^2}{2}\right\}H_{2n+1}\left(\frac{\eta}{\sqrt{2}}\right)
\exp\left\{-mi\omega\int_0^s\frac{ds'}{B(s')}+
2n\int_0^s\frac{\sigma(s')\xi^2(s')}{B(s')}\,ds'\right\}.
\end{align}
Thus the expressions (\ref{lambdatau}), (\ref{tau}), and (\ref{moi}) define the spectrum as
\begin{align}
Sp(L)=\left\{\lambda_0(1+O(\eps)) ,\bigcup_{n\geq 0,
m=\pm1,\pm2,\ldots}\lambda_{m,n}(1+O(\eps))\right\}.
\end{align}
As in (\ref{u0}), the forward eigenfunctions $u_{n,m}(\y)$ are related to the backward eigenfunctions
$v_{n,m}(\y)=\tilde Q_{m,n}(\eta,s)$ by
\begin{align}
u_{n,m}(\y)\sim\exp\left\{-\frac{\psi(\y)}{\eps}\right\}\left[K_0(\y)+O(\sqrt{\eps})\right]\bar
v_{n,m}(\y),\label{uvnm}
\end{align}
where in the initial variable
\beqq
v_{n,m}(\y)&=&\ds\exp\left\{-\frac{[\rho(\y)\xi(s(\y))]^2}{2\eps}\right\}H_{2n+1}
\left(\frac{\rho(\y)\xi(s(\y))}{\sqrt{2\eps}}\right)\\
&\times& \exp\left\{-mi\omega\int\limits_0^{s(\y)}\frac{ds'}{B(s')}+ 2n\int\limits_0^{s(\y)}\frac{\sigma(s')\xi^2(s')}{B(s')}\,ds'\right\}.
\eeqq
which in the boundary layer coordinates has the form
\begin{align}
u_{n,m}(\y)\sim\exp\left\{-\frac{\eta^2}{2}\right\}\left[K_0(\y)+O(\sqrt{\eps})\right]\bar
v_{n,m}(\y).\label{uvnmbl}
\end{align}
With the proper normalization the eigenfunctions $\{u_{n,m}(\y)\}$ and $\{v_{n,m}(\y)\}$ form a bi-orthonormal system.
\section{Applications}\label{s:Apparent}
{The asymptotic theory of section \ref{s:HO-ev} applies to a well-known model in neurophysiology, proposed in \citep*{Holcman2006}.} {In the absence of sensory stimuli the cerebral cortex is continuously active}. An example of this spontaneous activity is the phenomenon of voltage
transitions between two distinct levels, called Up and Down states, observed simultaneously when
recoding from many neurons \citep*{Ferster}, \citep*{Cossart}. The mathematical model {proposed in \citep*{Holcman2006} for cortical dynamics that exhibits spontaneous transitions between Up- and Down- states is given by the stochastic dynamics
\beq
\dot x &=&\frac{1-x}{t_r}-U x(y-T)H(y-T),\nonumber\\
&&\label{HTmodel}\\
\dot y &=&-\frac{y}{\tau}+\frac{x U w_T}{\tau}(y-T)H(y-T)+\frac{\sigma}{\sqrt{\tau}}\,\dot
w,\nonumber
\eeq
}where $x$ is a dimensionless synaptic depression parameter, $y$ is the membrane voltage, $U$ and
$t_r$ are utilization parameter and recovery time constant, respectively, $w_T$ is synaptic
strength, $\tau$ is a voltage time scale, $\sigma$ is noise amplitude, $H(\cdot)$ is the Heaviside
unit step function, and $\dot w$ is standard Gaussian white noise. {The model (\ref{HTmodel}) predicts that in a certain range of parameters the noiseless dynamics (when $\sigma=0$) has two basins of attractions: one around a focus, {which corresponds to an Up-state, and the second one is that of a stable equilibrium state, which corresponds to a Down-state. The basins of attraction are separated by an unstable limit cycle.}

{Figure \ref{f:Focus-Holcman-Tsodyks}A shows trajectories of (\ref{HTmodel}) that rotate several times around the focus before exiting the domain of attraction of the focus. Figure \ref{f:Focus-Holcman-Tsodyks}B shows the histogram of exit times oscillates with multiple peaks, as predicted by the theory presented in section \ref{s:HO-ev} above. The approximation of the histogram of exit times (\ref{pdf}) by the sum of the first two exponentials,
\beq
f(t)=A \exp (-\lambda_0 t) +B\exp (-\lambda_1t)\cos (\omega t+\phi),\label{f(t)}
\eeq
where $\lambda_0=\bar\tau^{-1}=1/2.46$ (computed empirically), $\lambda_1=2.5$ and the frequency is that of the focus (the
imaginary part of the Jacobian at the focus) $\omega=10.4$. The other parameters are $A=470$,
$B=600$, and $\phi=1.4$ (obtained by a numerical fit). The approximation (\ref{f(t)}) captures the first three oscillations that are smeared out in the exponentially decaying tail. The construction
of the short-time histogram requires the entire series expansion in (\ref{pdf}).}
\begin{figure}[http!]
\includegraphics[width=0.6\textwidth]{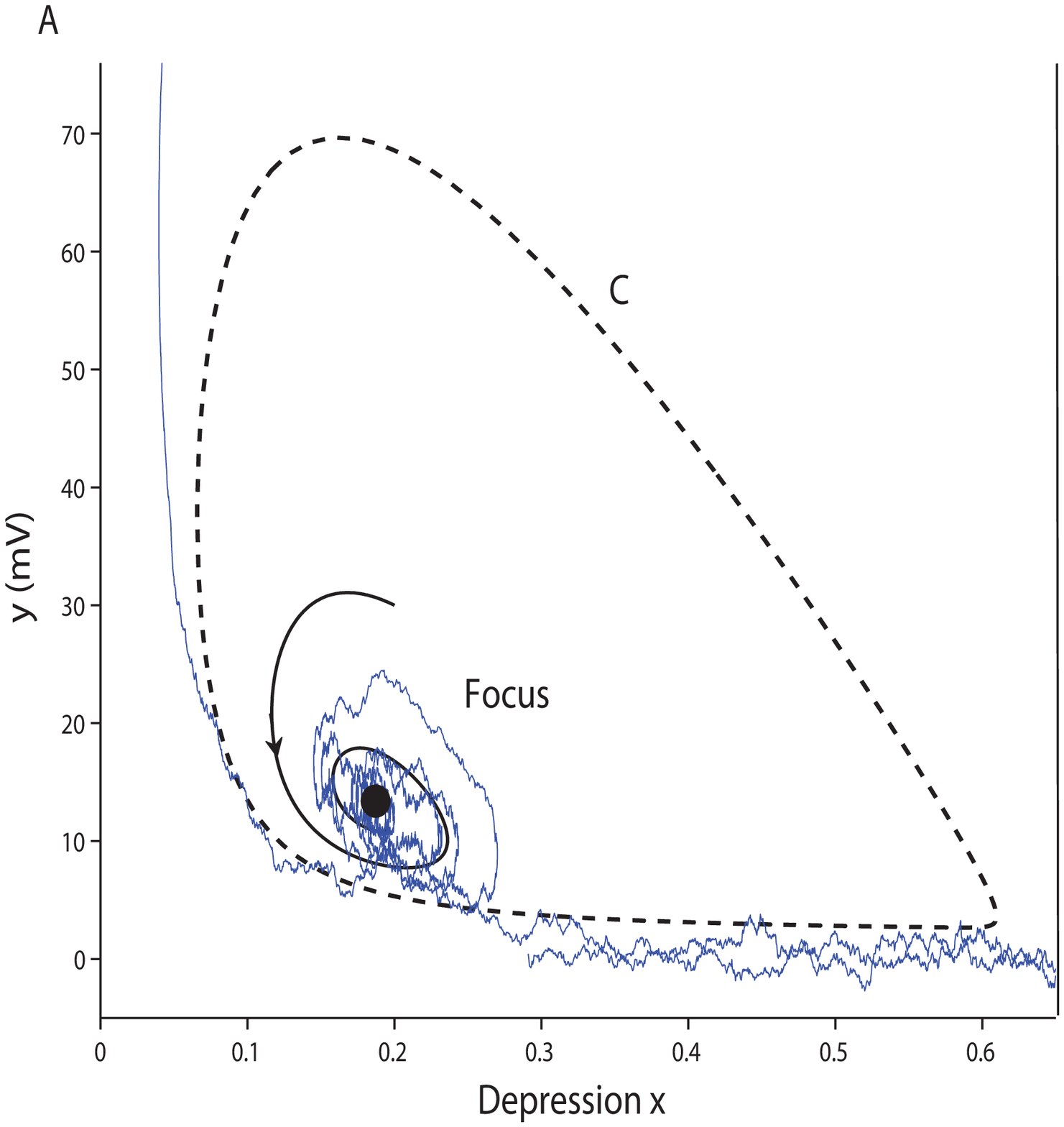}
\includegraphics[width=0.6\textwidth]{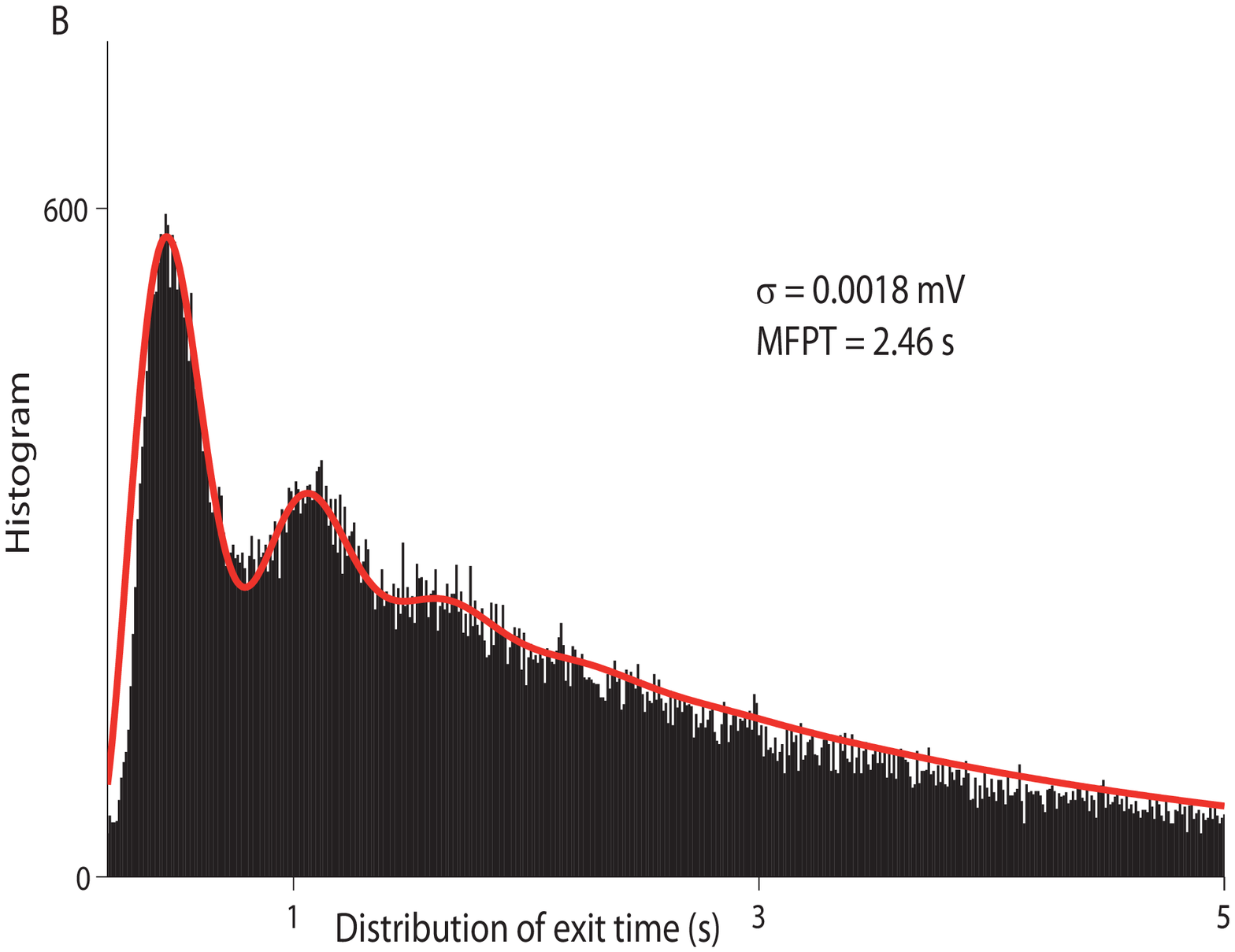}
\caption{\small { \bf The phase-plane dynamics of (\ref{HTmodel}), restricted to the Up-state.} {\bf A:} The unstable limit cycle is C (dashed line) and simulated trajectories (blue). The parameters are $\tau=0.05\,sec,t_r=0.8\,sec,U=0.5,w_T=12.6\,mV/Hz,T=2.0\,mV$.  {\bf B:} Histogram of exit times and its approximation by the first two terms of the expansion (\ref{pdf}) (marked red).} \label{f:Focus-Holcman-Tsodyks}
\end{figure}
{As indicated in section \ref{s:HO-ev}, the oscillation in the pdf} of exit times is a manifestation of the complex eigenvalues of the non-self adjoint Dirichlet problem for the corresponding Fokker-Planck operator inside the limit cycle.
\section{Summary and discussion}
This paper explains the oscillatory decay of the survival probability of the stochastic dynamics (\ref{SDEo}) that is activated over the boundary of the domain of attraction $D$ of the stable focus of the drift $\mb{a}(\x)$ by the small noise $\sqrt{2\eps}\,\mb{b}(\x_\eps(t))\, \dot\w(t)$. The boundary $\p D$ of the domain is an unstable limit cycle of $\mb{a}(\x)$. It is shown that the oscillations are not due a mysterious synchronization, but rather to complex eigenvalues of the Dirichlet problem for the Fokker-Planck operator in $D$. These are evaluated by a singular perturbation expansion of the spectrum of the non-self adjoint operator. The exact formula for the eigenvalues comes from the local expansion of the boundary layer in the neighborhood of the limit cycle. The expansion of the eigenvalues identifies for the first time the full and explicit spectrum of a non-self adjoint elliptic boundary value problem.

Oscillatory decay is manifested experimentally in the appearance of Up and Down states in the spontaneous activity of the cerebral cortex and in the simulations of its mathematical models \citep*{Holcman2006}. The oscillations are due to the competition between the driving noise and the underlying dynamical system.


\end{document}